\documentclass[12pt]{amsart}
\usepackage{graphicx,amsmath,amssymb}

\numberwithin{equation}{section}

\newcommand{\C}{\mathbb C}
\newcommand{\R}{\mathbb R}
\newcommand{\Z}{\mathbb Z}
\newcommand{\N}{\mathbb N}
\renewcommand{\d}{\prime} 
\newcommand{\dd}{{\prime \prime}}

\renewcommand{\Re}{{\rm Re}\,}
\renewcommand{\Im}{{\rm Im}\,}
\newtheorem{theorem}{Theorem}[section]
\newtheorem{lemma}[theorem]{Lemma}
 
\newtheorem{corollary}[theorem]{Corollary}

\newtheorem{definition}{Definition}

\newtheorem*{remarks}{Remarks}

\begin{document}
\title[]
{Half-Line non-self-adjoint Schr\"odinger operators  with polynomial potentials: Asymptotics of  eigenvalues}
\author[]
{Kwang C. Shin}
\address{Department of Mathematics, University of Missouri, Columbia, MO 65211, USA}
\date{February 16, 2005}

\begin{abstract}
For integers $m\geq 3$, we study the non-self-adjoint eigenvalue problems $-u^\dd(x)+(x^m+P(x))u(x)=E u(x)$, $0\leq x<+\infty$, with the boundary conditions  $u(+\infty)=0$ and $\alpha u(0)+\beta u^\d(0)=0$ for some $\alpha,\,\beta\in\C$ with $|\alpha|+|\beta|\not=0$, where $P(x)=a_1 x^{m-1}+a_2 x^{m-2}+\cdots+a_{m-1} x$ is a polynomial. We provide asymptotic expansions of the eigenvalue counting function and the eigenvalues $E_{n}$. Then we apply these to the inverse spectral problem, reconstructing some coefficients of polynomial potentials from asymptotic expansions of the eigenvalues.
\end{abstract}

\maketitle

\begin{center}
{\it Preprint.}
\end{center}

\baselineskip = 18pt

\section{Introduction}
\label{introduction}
In this paper, we study non-self-adjoint  Schr\"odinger operators in $L^2([0,+\infty))$, with monic  polynomial potentials of degree $m\geq 3$  and provide explicit asymptotic expansions of the eigenvalue counting functions and the eigenvalues $E_n$. Conversely, we reconstruct some coefficients of polynomial potentials from asymptotic expansions of the eigenvalues.

For an integer $m\geq 3$  and $(\alpha, \beta)\in\C^2\setminus\{(0,0)\}$, we consider the non-self-adjoint eigenvalue problems
\begin{align}\label{ptsym}
\left(H_{P}^{\alpha,\beta} u\right)(x):=\left[-\frac{d^2}{dx^2}+x^m+P(x)\right]u(x)=E u(x),\quad 0\leq x<+\infty,
\end{align}
for some $E\in\C$, with the boundary condition
\begin{equation}\label{bdcond}
\alpha u(0)+\beta u^\d(0)=0\quad\text{and \,\, $u(+\infty)=0$,}
\end{equation}
where $P$ is a polynomial of degree at most $m-1$ of the form 
\begin{equation}\nonumber
P(x)=a_1x^{m-1}+a_2x^{m-2}+\cdots+a_{m-1}x,\quad a_j\in\C\,\,\text{\,for $1\leq j\leq m-1$}.
\end{equation}

If a nonconstant function $u$ satisfies \eqref{ptsym} with some $E\in\C$ and the boundary condition \eqref{bdcond}, then we call $E$ an {\it eigenvalue} of $H_{P}^{\alpha,\beta}$ and $u$ an {\it eigenfunction of $H_{P}^{\alpha,\beta}$ associated with the eigenvalue $E$}. Also, the {\it geometric multiplicity of an eigenvalue $E$} is the number of linearly independent eigenfunctions  associated with the eigenvalue $E$. 

We number the eigenvalues $\{E_{n}\}_{n\geq n_0}$ in the order of nondecreasing magnitudes, counting their ``algebraic multiplicities'', where the integer $n_0$ could  depend on the potential and the boundary condition. In Theorem \ref{main_thm1} we show that for every large $n\in\N$, there exists $E_n$ satisfying  \eqref{main_result} below. However, we do not know the number of eigenvalues ``near'' zero, and this is why we need the number $n_0$.

Throughout this paper, we  use $E_n$ to denote  the eigenvalues $E_n=E_n(m,P,\alpha,\beta)$ of $H_{P}^{\alpha,\beta}$, without explicitly indicating their dependence on the potential and the boundary condition. 
Also, we let $$a:=(a_1,a_2,\ldots, a_{m-1})\in\C^{m-1}$$ be the coefficient vector of $P$. 

Before we state our main theorems, we first introduce some known facts by Sibuya \cite{Sibuya} about the eigenvalues $E_n$ of $H_{P}^{\alpha,\beta}$.
\begin{theorem}\label{main2}
The eigenvalues $E_{n}$ of $H_{P}^{\alpha,\beta}$ have the following properties.
\begin{enumerate}
\item[(I)] The set of all eigenvalues is a  discrete set in $\C$.
\item[(II)] The geometric multiplicity of every eigenvalue is one.
\item[(III)] Infinitely many eigenvalues, accumulating at infinity, exist.
\end{enumerate} 
\end{theorem}

This paper contains results on direct and inverse spectral problems. Theorem \ref{main_thm1} below is the main result, regarding asymptotic expansions of  ``eigenvalue counting functions''. The other results stated below in the Introduction are  deduced from   Theorem \ref{main_thm1}.

\subsection*{Direct spectral problem}
Here,
we first introduce the following theorem, regarding asymptotic expansions of a kind of eigenvalue counting functions, where we use multi-index notations with
$$\xi=(\xi_1,\xi_2,\dots,\xi_{m-1})\in\left(\N\cup\{0\}\right)^{m-1},\quad \text{and}\,\,\, \eta=(1,2,\dots,{m-1}).$$  Also, we use $|\xi|=\xi_1+\xi_2+\cdots+\xi_{m-1}$, $\xi!= \xi_1!\xi_2!\cdots\xi_{m-1}!$ and $a^{\xi}=a_1^{\xi_1}a_2^{\xi_2}\cdots a_{m-1}^{\xi_{m-1}}$. Also, $\lfloor x\rfloor$ is the largest integer that is less than or equal to $x\in\R$.
\begin{theorem}\label{main_thm1}
For  $a\in\C^{m-1}$, the eigenvalues $E_n$ of $H_{P}^{\alpha,\beta}$ satisfy
\begin{equation}\label{main_result}
\frac{1}{\pi}\sum_{j=0}^{\lfloor\frac{m+2}{2}\rfloor}d_j(a)E_n^{\frac{1}{2}+\frac{1-j}{m}}+o(1) =\left\{
\begin{array}{cl}n-\frac{1}{4},
 &\,\,\text{ if $\beta=0$,}\\
&\\
n+\frac{1}{4},
 &\,\,\text{ if $\beta\not=0$,}
\end{array}\right.
\end{equation}
 as $n\to+\infty$, where
 the error  term is uniform on any compact set of $a\in\C^{m-1}$ and 
\begin{equation}
d_{j}(a)=\left\{\begin{array}{cl}
\cos\left(\frac{(j-1)\pi}{m}\right)K_{m,j}(a)&\,\,\text{ if $0\leq j\leq\frac{m+1}{2}$,}\\
&\\
-\frac{\nu(a)}{m}\pi&\,\,\text{ if $m$ is even and $j=\frac{m+2}{2}$},
\end{array}\right.
\end{equation}
where
\begin{equation}\label{K_deF}
K_{m,0}(a)=K_{m,0,0}=\frac{B\left(\frac{1}{2},\,1+\frac{1}{m}\right)}{2\cos\left(\frac{\pi}{m}\right)},\,\,\,\, K_{m,j}(a)=\sum_{k=1}^jb_{j,k}(a)K_{m,j,k},\,\,1\leq j\leq \frac{m+2}{2}.
\end{equation}
Here $B(\cdot,\cdot)$ is the beta function and
\begin{equation}\label{K_definition}
K_{m,j,k}=\left\{\begin{array}{rl}
&\int_0^{\infty}\left(\frac{t^{mk-j}}{\left(t^m+1\right)^{k-\frac{1}{2}}}-t^{\frac{m}{2}-j}\right)\,dt,\,\,\text{ if }\,1\leq k\leq j\leq \frac{m+1}{2}\,\,\text{or}\,\,k=j=0,\\
&\int_0^{\infty}\left(\frac{t^{mk-\frac{m}{2}-1}}{\left(t^m+1\right)^{k-\frac{1}{2}}}-\frac{1}{t+1}\right)\,dt,\,\,\text{ if $m$ is even and }\,1\leq k\leq j= \frac{m+2}{2},
\end{array}
\right.
\end{equation} 
\begin{equation}\label{bjk_def}
b_{j,k}(a)={\frac{1}{2}\choose{k}}\sum_{\substack{|\xi|=k\\ \xi\cdot\eta=j}}\frac{k!}{\xi!}\,a^{\xi},\quad 1\leq k\leq j\leq \frac{m+2}{2},
\end{equation}
\begin{equation}\nonumber
\nu(a)=\left\{
                    \begin{array}{cl}
                  \sum_{k=1}^{\frac{m}{2}+1}b_{\frac{m}{2}+1,k}(a) \quad &\text{if $m$ is even,}\\
&\\
                   0\quad &\text{if $m$ is odd.}
                    \end{array}\right.
\end{equation} 
\end{theorem}
One can compute $K_{m,j,k}$ directly (or see \cite{Shin3}):
\begin{equation}\nonumber
K_{m,j,k}
=\left\{
                    \begin{array}{cl}
-\frac{2}{m}
\quad &\text{if $j=k=1$},\\
&\\
                  -\frac{2k-1}{m+2-2j}B\left(k-\frac{j-1}{m},\,\frac{1}{2}+\frac{j-1}{m}\right)   \quad &\text{if $1\leq k\leq j\leq\frac{m+1}{2}$, $j\not=1$},\\
&\\
                 \frac{2}{m}\left(\ln 2-\frac{1}{1}-\frac{1}{3}-\dots-\frac{1}{2k-5}-\frac{1}{2k-3}\right)  \quad &\text{if $m$ is even, $1\leq k\leq j=\frac{m+2}{2}$.}
                    \end{array}\right.\nonumber
\end{equation}

We obtain  \eqref{main_result} by investigating the asymptotic expansions of an entire function (the Stokes multiplier) whose zeros are the eigenvalues. In this paper, the ``algebraic multiplicity'' of an eigenvalue is the multiplicity of the zero of the Stokes multiplier.

Next, we let $N(t)$, $t\in\R$, be the eigenvalue counting function, that is, $N(t)$ is the number of eigenvalues $E$ of $H_{P}^{\alpha,\beta}$ such that $|E|\leq t$.  Then the following theorem on an asymptotic expansion of the eigenvalue counting function is a consequence of Theorem \ref{main_thm1}.
\begin{theorem}\label{main_thm2}
Let $a\in\C^{m-1}$ be fixed. Suppose that $\Im\left(K_{m,j}(a)\right)=0$ for $1\leq j\leq\frac{m+2}{2}$. Then $N(t)$ has the asymptotic expansion 
\begin{equation}\label{count_ft}
N(t)=\frac{1}{\pi}\sum_{j=0}^{\lfloor\frac{m+1}{2}\rfloor}\cos\left(\frac{(j-1)\pi}{m}\right)K_{m,j}(a)t^{\frac{1}{2}-\frac{j-1}{m}}+O(1),\quad\text{as $t\to+\infty$,}
\end{equation}
 where the error $O(1)$ is uniform for any compact set of $a\in\C^{m-1}$. 
\end{theorem}
\begin{proof}
In Corollary \ref{ineq_eq}, we show that $|E_n|<|E_{n+1}|$ for all large $n\in\N$.

Suppose that $|E_n|\leq t<|E_{n+1}|$. Then since for $s\in\R$,
$$\left(n+1\pm\frac{1}{4}\right)^{s}=\left(n\pm\frac{1}{4}\right)^{s}+O\left(n^{s-1}\right), \quad\text{as $n\to\infty$},$$
we see from Theorem \ref{eigen_asy} below that $|E_{n+1}|-|E_{n}|=O\left(n^{\frac{m-2}{m+2}}\right)$.
Thus, 
$$E_n^{\frac{1}{2}-\frac{j-1}{m}}=t^{\frac{1}{2}-\frac{j-1}{m}}\left(1-\frac{t-E_n}{t}\right)^{\frac{1}{2}-\frac{j-1}{m}}=t^{\frac{1}{2}-\frac{j-1}{m}}\left(1+O\left(\frac{t-E_n}{t}\right)\right)=t^{\frac{1}{2}-\frac{j-1}{m}}+O\left(1\right).$$
Hence, replacing $E_n^{\frac{1}{2}-\frac{j-1}{m}}$ in \eqref{main_result} by $t^{\frac{1}{2}-\frac{j-1}{m}}+O\left(1\right)$ and  solving the resulting equation for $n$ complete the proof.
\end{proof}

Next, from \eqref{main_result} we get $E_n$ in terms of $n$. 
\begin{theorem}\label{eigen_asy}
For each $a\in\C^{m-1}$, there exist some constants $e_{j}(a)\in\C$, $2\leq j\leq \frac{m+2}{2}$, such that
\begin{equation}\label{main_asyeq}
E_{n}=E_{n,0}+\sum_{j=1}^{\lfloor\frac{m+2}{2}\rfloor}e_{j}(a)E_{n,0}^{1-\frac{j}{m}}+o\left(E_{n,0}^{1-\frac{1}{m}\lfloor\frac{m+2}{2}\rfloor}\right),\quad\text{as $n\to+\infty$},
\end{equation}
where the error term is uniform for any compact set of $a\in\C^{m-1}$ and where
\begin{equation}\nonumber
E_{n,0}=\left(\frac{2\sqrt{\pi}\Gamma\left(\frac{3}{2}+\frac{1}{m}\right)}{\Gamma\left(1+\frac{1}{m}\right)}\right)^{\frac{2m}{m+2}}\times\left\{\begin{array}{lr}
\left(n-\frac{1}{4}\right)^{\frac{2m}{m+2}},\quad\text{if $\beta=0$},\\
&\\
\left(n+\frac{1}{4}\right)^{\frac{2m}{m+2}},\quad\text{if $\beta\not=0$},
\end{array}\right.
\end{equation}
and $e_j(a)$, $0\leq j\leq \frac{m+2}{2}$, are defined recurrently by $e_0(a)=1$ and
\begin{align}
&e_{j}(a)=-\frac{2m}{m+2}\left(\frac{d_j(a)}{d_0(a)}+\sum_{\substack{|\xi|=k\geq 2\\ \xi\cdot\eta=j}}{\frac{1}{2}+\frac{1}{m}\choose k}\frac{k!}{\xi!}\, e(a)^{\xi} +\sum_{r=1}^{j-1}\frac{d_r(a)}{d_0(a)}\sum_{\substack{|\xi|=k\\ \xi\cdot\eta=j-r}}{\frac{1}{2}+\frac{1-r}{m}\choose k}\frac{k!}{\xi!}\, e(a)^{\xi}\right),\nonumber
\end{align}
where $e(a)=(e_{1}(a),e_{2}(a),\dots,e_{m-1}(a))$.
\end{theorem}
We note, for the first summation  in the definition of $e_{j}(a)$ above, that $\xi\cdot\eta=j$ implies $\xi_{\ell}=0$ whenever $\ell\geq j$. Also, for the second summation, we point out that $\xi\cdot\eta=j-r\leq j-1$ implies $\xi_{\ell}=0$ whenever $\ell\geq j$.

When $P$ is real (i.e., $a\in\R^{m-1}$) and $x^m+P(x)$ is increasing and convex downwards on $[0,+\infty)$, Titchmarsh \cite[Chap.\ 7]{TIT} showed that
\begin{equation}\label{titcount}
N(t)\underset{t\to\infty}{=}\frac{1}{\pi}\int_0^{x_0}\sqrt{t-x^m-P(x)}\,dx+O(1),
\end{equation}
where $x_0=x_0(t)>0$ such that $t=x_0^m+P(x_0)$, provided that $\alpha=0$ or $\beta/\alpha$ real. Then from \eqref{titcount} one could get \eqref{count_ft}
  and hence  \eqref{main_result}. 

Voros \cite{VOR2} (cf.\ \cite{VOR3}) studied \eqref{ptsym} with arbitrary real polynomials $P$ under Dirichlet ($\beta=0$) and Neumann ($\alpha=0$) boundary conditions at $x=0$, and computed $d_0(a)$ and $d_1(a)$ explicitly. 

Fedoryuk \cite[\S 3.3]{Fedoryuk} considered \eqref{ptsym} with complex polynomial potentials and showed the existence of asymptotic expansions of the eigenvalues to all orders. Also, he computed $E_{n,0}$ explicitly.
 However, to the best of my knowledge Theorem \ref{main_thm1} in this generality does not appear in the literature to the date. 

Regarding monotonicity of modulus of $E_n$ for all large $n\in\N$.
\begin{corollary}\label{ineq_eq}
For each $a\in\C^{m-1}$ there exists $M>0$ such that $|E_{n}|<|E_{n+1}|$ if $n\geq M$. 
\end{corollary}
\begin{proof}
This is a consequence of  Theorem \ref{eigen_asy}. Or one can see that proof of Theorem 3 in \cite{Shin2} can be easily adapted for this case.
\end{proof}
\subsection*{Inverse spectral problem}
Here, we introduce results on inverse spectral problems, but first the following corollary is an easy consequence of Theorems \ref{main_thm1} and \ref{eigen_asy}, regarding how the coefficients of the asymptotic expansions depend on $a\in\C^{m-1}$.
\begin{corollary}\label{cor_poly}
Let $1\leq j\leq \frac{m+2}{2}$ be  a fixed integer. Then 
\begin{enumerate}
\item[(i)] $d_{j}(a)$ and $e_{j}(a)$ are  polynomials in $a_1,a_2,\dots,a_{j-1}, a_j$. In particular,  $d_{j}(a)$ and $e_{j}(a)$ are nonconstant linear functions of  $a_j$.
\item[(ii)]  $d_{j}(a)$ and $e_{j}(a)$ do not depend on $a_{j+1},a_{j+2},\dots, a_{m-1}$.
\end{enumerate}
\end{corollary}
\begin{proof}
Statements on $d_{j}(a)$ are direct consequences of the definition of $d_{j}(a)$ in Theorem \ref{main_thm1}. One can use statements on $d_{j}(a)$ and induction on $j$ to prove  statements on $e_{j}(a)$.
\end{proof}
Next, one can reconstruct some coefficients of the polynomial potential from the asymptotic expansion of the eigenvalues.
\begin{theorem}\label{inverse_thm}
 Let $1\leq j\leq\frac{m+1}{2}$ be a fixed integer. Then the asymptotic expansions of the eigenvalues $E_{n}$ of $H_{P}^{\alpha,\beta}$ of type \eqref{main_asyeq} with an error term $o\left(n^{\frac{2m-2j}{m+2}}\right)$ uniquely and explicitly determine  $a_k$ for all $1\leq k\leq j$. 
\end{theorem}
\begin{proof}
From the asymptotic expansion of the eigenvalues, one gets  $e_k(a)$ as an explicit polynomial in $a_1, a_2,\dots,a_k$ for every $1\leq k\leq j$.
Then since $e_{k}(a)$ is a nonconstant linear function of  $a_k$ and since $e_k(a)$ does not depend on $a_{\ell}$, $\ell>k$,  all $a_1, a_2,\dots,a_j$ can be found uniquely and explicitly.
\end{proof}
When $m$ is even, $j=\frac{m+2}{2}$ is allowed in Corollary \ref{cor_poly} while it is not allowed in Theorem \ref{inverse_thm}. This is due to the fact that our method in this paper does not determine the number $n_0$ in $\{E_{n}\}_{n\geq n_0}$.

\section{Properties of the solutions}
\label{prop_sect}
In this section, we introduce work of Hille \cite{Hille} and Sibuya \cite{Sibuya} about properties of the solutions of \eqref{ptsym}.

\begin{figure}[t]
    \begin{center}
    \includegraphics[width=.4\textwidth]{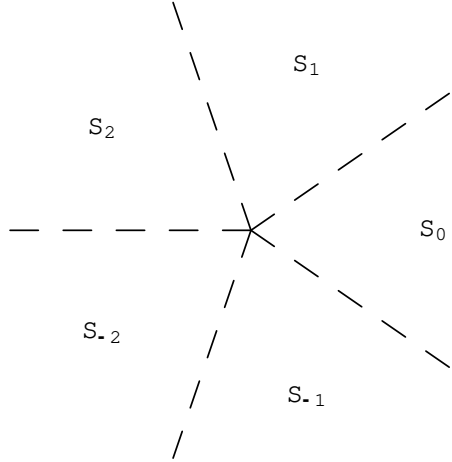}
    \end{center}
 \vspace{-.5cm}
\caption{The Stokes sectors for $m=3$. The dashed rays represent $\arg z=\pm\frac{\pi}{5},\,\pm\frac{3\pi}{5},\, \pi.$}\label{f:graph1}
\end{figure}

We first set 
$$\lambda=-E$$
and extend \eqref{ptsym} to the complex plane so that if $u$ is a solution of \eqref{ptsym} then
\begin{equation}\label{rotated}
-u^\dd(z)+[z^m+P(z)+\lambda]u(z)=0,\quad z\in\C.
\end{equation} 
It is known that solutions of \eqref{rotated} have rather simple asymptotic behavior near infinity in the complex plane \cite[\S 7.4]{Hille}. We will describe this simple asymptotic behavior of the solutions near infinity by using the following definition.
\begin{definition}
{\rm {\it The Stokes sectors} $S_k$ of the equation \eqref{rotated}
are
$$ S_k=\left\{z\in \C:\left|\arg (z)-\frac{2k\pi}{m+2}\right|<\frac{\pi}{m+2}\right\}\quad\text{for}\quad k\in \Z.$$ }
\end{definition}
See Figure \ref{f:graph1}.

Hille \cite[\S 7.4]{Hille} showed that every nonconstant solution of \eqref{rotated} either decays to zero or blows up exponentially, in each Stokes sector $S_k$. 
\begin{lemma}[\protect{\cite[\S 7.4]{Hille}}]\label{gen_pro}
${} $

\begin{itemize}
\item [(i)] For each $k\in\Z$, every solution $u$ of \eqref{rotated}  is asymptotic to 
\begin{equation} \label{asymp-formula}
(const.)z^{-\frac{m}{4}}\exp\left[\pm \int^z \left[z^m+P(z)+\lambda\right]^{\frac{1}{2}}\,dz\right] 
\end{equation}
as $z \rightarrow \infty$ in every closed subsector of $S_k$.

\item [(ii)] If a nonconstant solution $u$ of \eqref{rotated} decays in $S_k$, it must blow up in $S_{k-1}\cup S_{k+1}$. However, when $u$ blows up in $S_k$, $u$ need not be decaying in $S_{k-1}$ or in $S_{k+1}$.
\end{itemize}
\end{lemma}
Lemma \ref{gen_pro} (i) implies that if $u$ decays along one ray in $S_k$, then it decays along all rays in $S_k$. Also, if $u$ blows up along one ray in $S_k$, then it blows up along all rays in $S_k$. 

We will use
$$\omega=\exp\left[\frac{2\pi i}{m+2}\right]$$
and  we define  
\begin{equation}\nonumber
b_j(a)=\sum_{k=1}^j b_{j,k}(a),\quad 1\leq j\leq \frac{m+2}{2}.
\end{equation}
 We further define $r_m=-\frac{m}{4}$ if $m$ is odd, and $r_m=-\frac{m}{4}-b_{\frac{m}{2}+1}(a)$ if $m$ is even. 

Now we are ready to introduce some results of Sibuya \cite{Sibuya} that is the main ingredient of the proof of Theorem \ref{main_thm1}. 
\begin{theorem}\label{prop}
Equation \eqref{rotated}, with $a\in \C^{m-1}$, admits a solution  $f(z,a,\lambda)$ with the following properties.
\begin{enumerate}
\item[(i)] $f(z,a,\lambda)$ is an entire function of $z,a $ and $\lambda$.
\item[(ii)] $f(z,a,\lambda)$ and $f^\d(z,a,\lambda)=\frac{\partial}{\partial z}f(z,a,\lambda)$ admit the following asymptotic expansions. Let $\varepsilon>0$. Then
\begin{align}
f(z,a,\lambda)=&\qquad z^{r_m}(1+O(z^{-1/2}))\exp\left[-F(z,a,\lambda)\right],\nonumber\\
f^\d(z,a,\lambda)=&-z^{r_m+\frac{m}{2}}(1+O(z^{-1/2}))\exp\left[-F(z,a,\lambda) \right],\nonumber
\end{align}
as $z$ tends to infinity in  the sector $|\arg z|\leq \frac{3\pi}{m+2}-\varepsilon$, uniformly on each compact set of $(a,\lambda)$-values . 
Here
\begin{equation}\nonumber
F(z,a,\lambda)=\frac{2}{m+2}z^{\frac{m}{2}+1}+\sum_{1\leq j<\frac{m}{2}+1}\frac{2}{m+2-2j}b_j(a) z^{\frac{1}{2}(m+2-2j)}.
\end{equation}
\item[(iii)] For each fixed $a\in\C^{m-1}$ and $\delta>0$, $f$ and $f^\d$ also admit the asymptotic expansions, 
\begin{align}
f(0,a,\lambda)=&[1+o(1)]\lambda^{-1/4}\exp\left[L(a,\lambda)\right],\label{eq1}\\
f^\d(0,a,\lambda)=&-[1+o(1)]\lambda^{1/4}\exp\left[L(a,\lambda)\right],\label{eq2}
\end{align}
as $\lambda\to\infty$ in the sector $|\arg(\lambda)|\leq\pi-\delta$, uniformly on each compact set of $a\in\C^{m-1}$, where 
\begin{align}
L(a,\lambda)=\left\{
                    \begin{array}{rl}
                    &\int_0^{+\infty}\left(\sqrt{t^m+P(t)+\lambda}- t^{\frac{m}{2}}-\sum_{j=1}^{\frac{m+1}{2}}b_j(a)t^{\frac{m}{2}-j}\right)\,dt \quad \text{if $m$ is odd,}\\
                   &\int_0^{+\infty}\left(\sqrt{t^m+P(t)+\lambda}- t^{\frac{m}{2}}-\sum_{j=1}^{\frac{m}{2}}b_j(a)t^{\frac{m}{2}-j}-\frac{b_{\frac{m}{2}+1}(a)}{t+1}\right)\,dt  \quad \text{if $m$ is even.}
                    \end{array}
                         \right. \nonumber
\end{align}
\item[(iv)] The entire functions  $\lambda\mapsto f(0,a,\lambda)$ and $\lambda\mapsto f^\d(0,a,\lambda)$ have orders $\frac{1}{2}+\frac{1}{m}$.
\end{enumerate}
\end{theorem}
\begin{proof}
In  Sibuya's book \cite{Sibuya}, see Theorem 6.1 for a proof of (i) and (ii), and Theorem 19.1 for a proof of (iii).  Moreover, (iv) is a consequence of (iii) along with Theorem 20.1 in \cite{Sibuya}. Note that properties (i), (ii), and (iii) are summarized on pages 112--113 of Sibuya \cite{Sibuya}. 
\end{proof}

\begin{remarks}
{\rm (I) Uniformity of the error term in Theorem \ref{main_thm1}  is essentially due to uniformity of error terms in \eqref{eq1} and \eqref{eq2}.

(II) In this paper  we will deal with numbers like $\left(\omega^{\nu}\lambda\right)^{s}$ for some $s\in\R$, and $\nu\in\C$. As usual, we will use
$$\omega^{\nu}=\exp\left[\nu \frac{2\pi i}{m+2}\right]$$
and if $\arg(\lambda)$ is specified, then
$$\arg\left(\left(\omega^{\nu}\lambda\right)^{s}\right)=s\left[\arg(\omega^{\nu})+\arg(\lambda)\right]=s\left[\Re(\nu)\frac{2\pi}{m+2}+\arg(\lambda)\right],\quad s\in\R.$$
If $s\not\in\Z$ then the branch of $\lambda^s$ is chosen to be the negative real axis.
}
\end{remarks}

In \cite{Shin2}, the following asymptotic expansion of $L(a,\cdot)$ is proved. 
\begin{lemma}\label{asy_lemma}
Let $m\geq 3$ and $a\in\C^{m-1}$ be fixed. Then 
\begin{equation}
L(a,\lambda)=\left\{\begin{array}{rl}
&\sum_{j=0}^{\frac{m+1}{2}}K_{m,j}(a)\lambda^{\frac{1}{2}+\frac{1-j}{m}}+O\left(|\lambda|^{-\frac{1}{2m}}\right)\,\,\text{if $m$ is odd,}\\
&\sum_{j=0}^{\frac{m}{2}+1}K_{m,j}(a)\lambda^{\frac{1}{2}+\frac{1-j}{m}}-\frac{b_{\frac{m}{2}+1}(a)}{m}\ln(\lambda)+O\left(|\lambda|^{-\frac{1}{m}}\right)\,\,\text{if $m$ is even,}
\end{array}
\right.\nonumber
\end{equation}
as $\lambda\to\infty$ in the sector $|\arg(\lambda)|\leq \pi-\delta$, uniformly on each compact set of $a\in\C^{m-1}$. 
\end{lemma}
\begin{proof}
See \cite{Shin2} for a proof.
\end{proof} 

Sibuya \cite{Sibuya} introduced solutions of \eqref{rotated} that decays in $S_k$, $k\in\Z$. Before we introduce this, we let 
\begin{equation}\label{G_def}
G^{\ell}(a):=(\omega^{-\ell}a_1, \omega^{-2\ell}a_2,\ldots,\omega^{-(m-1)\ell}a_{m-1})\quad \text{for}\quad \ell\in \Z.
\end{equation}
Then we have the following lemma, regarding properties of $G^{\ell}(\cdot)$.
\begin{lemma}\label{lemma_25}
For $a\in\C^{m-1}$ fixed, and $\ell_1,\ell_2,\ell\in\Z$,
$G^{\ell_1}(G^{\ell_2}(a))=G^{\ell_1+\ell_2}(a)$, and
\begin{equation}\nonumber
b_{j,k}(G^{\ell}(a))=\omega^{-j\ell}b_{j,k}(a),\quad \ell\in\Z.
\end{equation}
\end{lemma}

Next, recall that the function
 $f(z,a,\lambda)$  in Theorem \ref{prop} solves \eqref{rotated} and decays to zero exponentially as $z\rightarrow \infty$ in $S_0$, and  blows up in $S_{-1}\cup S_1$. One can check that the function
$$f_k(z,a,\lambda):=f(\omega^{-k}z,G^k(a),\omega^{2k}\lambda),\quad k\in\Z,$$
 which is obtained by scaling $f(z,G^k(a),\omega^{2k}\lambda)$ in the $z$-variable, also solves \eqref{rotated}. It is clear that $f_0(z,a,\lambda)=f(z,a,\lambda)$, and that $f_k(z,a,\lambda)$ decays in $S_k$ and blows up in $S_{k-1}\cup S_{k+1}$ since $f(z,G^k(a),\omega^{2k}\lambda)$ decays in $S_0$. Since no nonconstant solution decays in two consecutive Stokes sectors (see Lemma \ref{gen_pro} (ii)), $f_{0}$ and $f_{-1}$ are linearly independent and hence any solution of \eqref{rotated} can be expressed as a linear combination of these two. Especially,  there exist some coefficients $C(a,\lambda)$ and $\widetilde{C}(a,\lambda)$ such that
\begin{equation}\label{basic_eq}
f_{1}(z,a,\lambda)=C(a,\lambda)f_{0}(z,a,\lambda)+\widetilde{C}(a,\lambda)f_{-1}(z,a,\lambda).
\end{equation}

We then see that 
\begin{equation}\label{C_def}
C(a,\lambda)=\frac{W_{-1,1}(a,\lambda)}{W_{-1,0}(a,\lambda)}\quad\text{and}\quad \widetilde{C}(a,\lambda)=\frac{W_{1,0}(a,\lambda)}{W_{-1,0}(a,\lambda)},
\end{equation}
where $W_{j,\ell}=f_jf_{\ell}^\d -f_j^\d f_{\ell}$ is the Wronskian of $f_j$ and $f_{\ell}$. Since both $f_j,\,f_{\ell}$ are solutions of the same linear equation \eqref{rotated}, we know that the Wronskians are constant functions of $z$. Also, $f_k$ and $f_{k+1}$ are linearly independent, and hence $W_{k,k+1}\not=0$ for all $k\in \Z$. 

Moreover, we have the following lemma that is useful  later on.
\begin{lemma}\label{shift_lemma}
Suppose $k,\,j\in\Z$. Then
\begin{equation}\label{kplus1}
W_{k+1,j+1}(a,\lambda)=\omega^{-1}W_{k,j}(G(a),\omega^2\lambda),
\end{equation}
and $W_{0,1}(a,\lambda)=2\omega^{\mu(a)}$, where
\begin{eqnarray}
\mu(a)=\left\{
              \begin{array}{rl}
              \frac{m}{4}  \quad &\text{if $m$ is odd,}\\
              \frac{m}{4}- b_{\frac{m}{2}+1}(a) \quad &\text{if $m$ is even.}
              \end{array}
                         \right. \nonumber
\end{eqnarray}

\end{lemma}
\begin{proof}
See Sibuya \cite[pages 116-118]{Sibuya}.
\end{proof}
Thus, by Lemma \ref{shift_lemma},
\begin{equation}\nonumber
\widetilde{C}(a,\lambda)=\frac{W_{1,0}(a,\lambda)}{W_{-1,0}(a,\lambda)}=-\frac{2\omega^{\mu(a)}}{2\omega \omega^{\mu(G^{-1}(a))}}=-\omega^{-1-2\nu(a)},
\end{equation}
where $\nu(a)=\frac{m}{4}-\mu(a)$, that is, 
\begin{eqnarray}
\nu(a)=\left\{
                         \begin{array}{rl}
                         0 & \quad \text{if $m$ is odd,}\\
                         b_{\frac{m}{2}+1}(a)& \quad \text{if $m$ is even.}
                          \end{array}
                         \right. \label{def_nu}
\end{eqnarray}

\section{Asymptotics of $f(0,a,\lambda)$ and $f^\d(0,a,\lambda)$}\label{sec_3}
The asymptotics of $f(0,a,\lambda)$ and $f^\d(0,a,\lambda)$ as $\lambda\to\infty$ in the sector $|\arg(\lambda)|\leq \pi-\delta$ are given by \eqref{eq1} and \eqref{eq2}, respectively. In this section, we provide the  asymptotics of $f(0,a,\lambda)$ and $f^\d(0,a,\lambda)$ as $\lambda\to\infty$ in a  sector near the negative real axis.

In  \cite{Shin2}, we showed
 the following asymptotic expansion of $W_{-1,1}(a,\lambda)$ as $\lambda\to\infty$ in a sector near the negative real axis.
\begin{theorem}\label{thm_neg}
Let $m\geq 3$, $a\in\C^{m-1}$ and $0<\delta<\frac{\pi}{m+2}$ be fixed. Then
\begin{equation}\label{asy_1}
W_{-1,1}(a,\lambda)=[2i+o(1)]\exp\left[L(G^{-1}(a),\omega^{-2}\lambda)+L(G(a),\omega^{-m}\lambda)\right],
\end{equation}
as $\lambda\to \infty$ along the rays in the sector
\begin{equation}\label{sector1}
\pi-\frac{4\pi}{m+2}+\delta\leq \arg(\lambda)\leq \pi+\frac{4\pi}{m+2}-\delta, 
\end{equation}
where the error term is uniform on any compact set of $a\in\C^{m-1}$.
\end{theorem}
\begin{proof}
See \cite[Theorem 12]{Shin2} for a proof.
\end{proof}
We will use this in the next theorem, regarding asymptotics of $f(0,a,\lambda)$ and $f^\d(0,a,\lambda)$ near the negative real axis.
\begin{theorem}
Let $a\in\C^{m-1}$ be fixed. Then
\begin{align}
f(0,a,\lambda)=&\left(\frac{i}{2}\omega^{-\nu(a)}+o(1)\right)\lambda^{-\frac{1}{4}}\exp\left[-L(G^{-1}(a),\omega^{-2}\lambda)\right]\nonumber\\
&+\left(\frac{1}{2}\omega^{-3\nu(a)}+o(1)\right)\lambda^{-\frac{1}{4}}\exp\left[-L(G(a),\omega^{-m}\lambda)\right],\label{form_1}\\
f^\d(0,a,\lambda)=&\left(\frac{i}{2}\omega^{-\nu(a)}+o(1)\right)\lambda^{\frac{1}{4}}\exp\left[-L(G^{-1}(a),\omega^{-2}\lambda)\right]\nonumber\\
&-\left(\frac{1}{2}\omega^{-3\nu(a)}+o(1)\right)\lambda^{\frac{1}{4}}\exp\left[-L(G(a),\omega^{-m}\lambda)\right], \quad\text{as $\lambda\to\infty$ in  \eqref{sector1},}\label{form_2} 
\end{align}
where the error terms are uniform on any compact set of $a\in\C^{m-1}$ and where $\arg\left(\lambda^{\pm\frac{1}{4}}\right)=\pm\frac{1}{4}\arg(\lambda)$ in the sector \eqref{sector1}.
\end{theorem}
\begin{proof}
From \eqref{basic_eq} and \eqref{C_def}, and Lemma \ref{shift_lemma}, we have
\begin{align}
f(z,a,\lambda)&=f_{0}(z,a,\lambda)\nonumber\\
&=\frac{1}{C(a,\lambda)}\left[f_{1}(z,a,\lambda)-\widetilde{C}(a,\lambda)f_{-1}(z,a,\lambda)\right]\nonumber\\
&=\frac{2\omega^{1+\mu(G^{-1}(a)}}{W_{-1,1}(a,\lambda)}\left[f(\omega^{-1}z,G(a),\omega^2\lambda)+\omega^{-1-2\nu(a)}f(\omega z,G^{-1}(a),\omega^{-2}\lambda)\right].\label{trans_eq}
\end{align}
So we examine  asymptotics  of $f(\omega^{-1}z,G(a),\omega^2\lambda) +\omega^{-1-2\nu(a)}f(\omega z,G^{-1}(a),\omega^{-2}\lambda)$ and its derivative at $z=0$. Using \eqref{eq1} and the fact that 
$f(0,G(a),\omega^2\lambda)=f(0,G(a),\omega^{-m}\lambda)$, we have
\begin{align}
&f(0,G(a),\omega^2\lambda)+\omega^{-1-2\nu(a)}f(0,G^{-1}(a),\omega^{-2}\lambda)\nonumber\\
&=f(0,G(a),\omega^{-m}\lambda)+\omega^{-1-2\nu(a)}f(0,G^{-1}(a),\omega^{-2}\lambda)\nonumber\\
&=\left(1+o(1)\right)\left(\omega^{-m}\lambda\right)^{-\frac{1}{4}}\exp\left[L(G(a),\omega^{-m}\lambda)\right]\nonumber\\
&+\omega^{-1-2\nu(a)}\left(1+o(1)\right)\left(\omega^{-2}\lambda\right)^{-\frac{1}{4}}\exp\left[L(G^{-1}(a),\omega^{-2}\lambda)\right]\nonumber\\
&=\left(\omega^{\frac{m}{4}}+o(1)\right)\lambda^{-\frac{1}{4}}\exp\left[L(G(a),\omega^{-m}\lambda)\right]+\left(\omega^{-\frac{1}{2}-2\nu(a)}+o(1)\right)\lambda^{-\frac{1}{4}}\exp\left[L(G^{-1}(a),\omega^{-2}\lambda)\right],\nonumber
\end{align}
as $\lambda\to\infty$ in the sector \eqref{sector1}.
Then \eqref{form_1} is obtained from  \eqref{asy_1} and \eqref{trans_eq}.

Next, we differentiate \eqref{trans_eq} with respect to $z$ and evaluate the resulting equation at $z=0$ to get
\begin{equation}\label{trans_eq1}
f^\d(0,a,\lambda)=\frac{2\omega^{1+\mu(G^{-1}(a))}}{W_{-1,1}(a,\lambda)}\left[\omega^{-1}f^\d(0,G(a),\omega^{-m}\lambda)+\omega^{-2\nu(a)}f^\d(0,G^{-1}(a),\omega^{-2}\lambda)\right].
\end{equation}
Using \eqref{eq2}, we have
\begin{align}
&\omega^{-1}f^\d(0,G(a),\omega^{-m}\lambda)+\omega^{-2\nu(a)}f^\d(0,G^{-1}(a),\omega^{-2}\lambda)\nonumber\\
&=\left(i\omega^{-\frac{1}{2}}+o(1)\right)\lambda^{\frac{1}{4}}\exp\left[L(G(a),\omega^{-m}\lambda)\right]-\left(\omega^{-\frac{1}{2}-2\nu(a)}+o(1)\right)\lambda^{\frac{1}{4}}\exp\left[L(G^{-1}(a),\omega^{-2}\lambda)\right],\nonumber
\end{align}
as $\lambda\to\infty$ in the sector \eqref{sector1}.
Then  this along with \eqref{asy_1} and \eqref{trans_eq1} yields \eqref{form_2}.

Finally, the uniformity of the error terms in \eqref{trans_eq} and \eqref{trans_eq1} is due to the uniformity of the error terms in \eqref{eq1}, \eqref{eq2}, and \eqref{asy_1}.
\end{proof}
\section{Proof of Theorem \ref{main_thm1}}\label{sec_4}
In this section we prove Theorem \ref{main_thm1}.

\begin{proof}[Proof of Theorem ~\ref{main_thm1} for  Dirichlet boundary condition at $x=0$.]
From \eqref{form_1}
\begin{align}
&2\omega^{3\nu(a)}\lambda^{\frac{1}{4}}\exp\left[L(G^{-1}(a),\omega^{-2}\lambda)+o(1)\right]f(0,a,\lambda)\nonumber\\
&=\exp\left[L(G^{-1}(a),\omega^{-2}\lambda)-L(G(a),\omega^{-m}\lambda)+o(1)\right]
+i\omega^{2\nu(a)},\quad\text{as $\lambda\to\infty$ in  \eqref{sector1}.}\nonumber
\end{align}
Since 
\begin{align}
&L(G^{-1}(a),\omega^{-2}\lambda)-L(G(a),\omega^{-m}\lambda)\nonumber\\
&=K_{m}\left(\omega^{-2}\lambda\right)^{\frac{1}{2}+\frac{1}{m}}(1+o(1))-K_{m}\left(\omega^{-m}\lambda\right)^{\frac{1}{2}+\frac{1}{m}}(1+o(1))\nonumber\\
&=K_m\left(\exp\left[-\frac{2\pi}{m} i\right]-\exp\left[-\pi i\right]\right)\lambda^{\frac{1}{2}+\frac{1}{m}}(1+o(1))\nonumber\\
&=K_m\left(1+\exp\left[-\frac{2\pi}{m} i\right]\right)\lambda^{\frac{1}{2}+\frac{1}{m}}(1+o(1)),\nonumber
\end{align}
and since $\arg\left(1+\exp\left[-\frac{2\pi}{m} i\right]\right)=-\frac{\pi}{m}$, we have
\begin{equation}\nonumber
\arg\left(L(G^{-1}(a),\omega^{-2}\lambda)-L(G(a),\omega^{-m}\lambda)\right)=-\frac{\pi}{m}+\frac{m+2}{2m}\arg(\lambda)+o(1). 
\end{equation}

Thus, if $\pi-\frac{4\pi}{m+2}+\delta\leq \arg(\lambda)\leq \pi+\frac{4\pi}{m+2}-\delta$ and $|\lambda|$ is large, we have
\begin{equation}
\frac{\pi}{2}-\frac{2\pi}{m}+o(1)\leq \arg\left(L(G^{-1}(a),\omega^{-2}\lambda)-L(G(a),\omega^{-m}\lambda)\right)\leq\frac{\pi}{2}+\frac{2\pi}{m}+o(1).
\end{equation}
So $\lambda\mapsto L(G^{-1}(a),\omega^{-2}\lambda)-L(G(a),\omega^{-m}\lambda)$ maps the sector \eqref{sector1} near infinity onto a region containing  $|\arg(\lambda)-\frac{\pi}{2}|\leq \varepsilon_1$ and $|\lambda|\geq M_0$ for some positive real numbers $\varepsilon, M_0.$  Hence, there exists a sequence of the numbers $\lambda_n$ in \eqref{sector1} such that
\begin{equation}\label{bas_eqn}
L(G^{-1}(a),\omega^{-2}\lambda_n)-L(G(a),\omega^{-m}\lambda_n)+o(1)\underset{n\to+\infty}{=}2n\pi i+2\nu(a)\frac{2\pi i}{m+2}-\frac{\pi}{2}i,
\end{equation}
for all large $n\in\N$ so that $f(0,a,\lambda_n)=0$.
Next, by  \eqref{K_deF} and Lemma \ref{asy_lemma}
\begin{align}
&L(G^{-1}(a),\omega^{-2}\lambda_n)-L(G(a),\omega^{-m}\lambda_n)\nonumber\\
&=\sum_{j=0}^{\lfloor\frac{m}{2}+1\rfloor}\left(K_{m,j}(G^{-1}(a))(\omega^{-2}\lambda_n)^{\frac{1}{2}+\frac{1-j}{m}}-K_{m,j}(G(a))(\omega^{-m}\lambda_n)^{\frac{1}{2}+\frac{1-j}{m}}\right)\nonumber\\
&-\frac{\nu(G^{-1}(a))}{m}\ln(\omega^{-2}\lambda_n)+\frac{\nu(G(a))}{m}\ln(\omega^{-m}\lambda_n)+o(1)\nonumber
\end{align}
\begin{align}
\quad\quad\quad&=\sum_{j=0}^{\lfloor\frac{m}{2}+1\rfloor}\left(\omega^{j}\omega^{-2(\frac{1}{2}+\frac{1-j}{m})}-\omega^{-j}\omega^{-m(\frac{1}{2}+\frac{1-j}{m})}\right)K_{m,j}(a)\lambda_n^{\frac{1}{2}+\frac{1-j}{m}}\nonumber\\
&-\frac{\nu(a)}{m}\frac{4\pi}{m+2}i+\frac{\nu(a)}{m}\frac{2m\pi}{m+2}i+o(1)\nonumber\\
&=2i\sum_{j=0}^{\lfloor\frac{m}{2}+1\rfloor}\sin\left(\frac{(m-2+2j)\pi}{2m}\right)K_{m,j}(a)(-\lambda_n)^{\frac{1}{2}+\frac{1-j}{m}}+\frac{\nu(a)}{m}\frac{(2m-4)\pi}{m+2}i+o(1).\nonumber
\end{align}
So this and \eqref{bas_eqn} yield
\begin{equation}\nonumber
2i\sum_{j=0}^{\lfloor\frac{m}{2}+1\rfloor}\sin\left(\frac{(m-2+2j)\pi}{2m}\right)K_{m,j}(a)(-\lambda_n)^{\frac{1}{2}+\frac{1-j}{m}}-\frac{2\nu(a)}{m}\pi i+o(1)\underset{n\to+\infty}{=}\left(2n-\frac{1}{2}\right)\pi i.
\end{equation}
Finally, we use $\sin(\pi/2+\theta)=\cos(\theta)$ and $E_n=-\lambda_n$ to complete the proof.
\end{proof}

Next, we prove Theorem ~\ref{main_thm1} for the case when $\beta\not=0$ in \eqref{bdcond}.
\begin{proof}[Proof of Theorem ~\ref{main_thm1} for other  boundary conditions]
Using \eqref{form_1} and \eqref{form_2}, one gets
\begin{align}
&\alpha f(0,a,\lambda)+\beta f^\d(0,a,\lambda)\nonumber\\
&=\left\{\frac{\alpha}{2\lambda^{\frac{1}{4}}}\left(i\omega^{-\nu(a)}+o(1)\right)+\frac{\beta\lambda^{\frac{1}{4}}}{2}\left(i\omega^{-\nu(a)}+o(1)\right)\right\}\exp\left[-L(G^{-1}(a),\omega^{-2}\lambda)\right]\nonumber\\
&+\left\{\frac{\alpha}{2\lambda^{\frac{1}{4}}}\left(\omega^{-3\nu(a)}+o(1)\right)-\frac{\beta\lambda^{\frac{1}{4}}}{2}\left(\omega^{-3\nu(a)}+o(1)\right)\right\}\exp\left[-L(G(a),\omega^{-m}\lambda)\right]\nonumber\\
&=\frac{\beta\lambda^{\frac{1}{4}}}{2}\Big\{\left(i\omega^{-\nu(a)}+o(1)\right)\exp\left[-L(G^{-1}(a),\omega^{-2}\lambda)\right]\label{form_3}\\
&\quad\quad\quad-\left(\omega^{-3\nu(a)}+o(1)\right)\exp\left[-L(G(a),\omega^{-m}\lambda)\right]\Big\},\nonumber
\end{align}
as $\lambda\to\infty$ in the sector \eqref{sector1}, where the error terms are uniform on any compact set of $a\in\C^{m-1}$ and where $\arg\left(\lambda^{\pm\frac{1}{4}}\right)=\pm\frac{1}{4}\arg(\lambda)$ in the sector \eqref{sector1}.
Since $\beta\not=0$,
\begin{align}
&\frac{2}{\beta\lambda^{\frac{1}{4}}}\exp\left[L(G^{-1}(a),\omega^{-2}\lambda)+o(1)\right]\left[\alpha f(0,a,\lambda)+\beta f^\d(0,a,\lambda)\right]\nonumber\\
&=i\omega^{-\nu(a)}-\omega^{-3\nu(a)}\exp\left[L(G^{-1}(a),\omega^{-2}\lambda)-L(G(a),\omega^{-m}\lambda)+o(1)\right],\nonumber
\end{align}
as $\lambda\to\infty$ in the sector \eqref{sector1}. 

Thus, like in the proof of Theorem \ref{main_thm1} for $\beta=0$, there exists a sequence of $\lambda_n$ such that 
\begin{equation}\nonumber
L(G^{-1}(a),\omega^{-2}\lambda_n)-L(G(a),\omega^{-m}\lambda_n)+o(1)\underset{n\to+\infty}{=}\left(2n+\frac{1}{2}\right)\pi i+2\nu(a)\frac{2\pi i}{m+2},
\end{equation}
for all large $n\in\N$ so that $\alpha f(0,a,\lambda_n)+\beta f^\d(0,a,\lambda_n)=0$. Here we have $\left(2n+\frac{1}{2}\right)\pi i$ in the place of $\left(2n-\frac{1}{2}\right)\pi i$ in \eqref{bas_eqn}. So one can complete the proof by following the methods in the proof for $\beta=0$ case.
\end{proof}
\section{Proof of Theorem \ref{eigen_asy}}\label{sec_9}
We will prove existence of $e_j(a)$ by induction on $j$. In doing so we will recurrently find  $e_j(a)$.

From \eqref{main_result}  we have
\begin{equation}\label{main_result1}
\sum_{j=0}^{\lfloor\frac{m}{2}+1\rfloor}\frac{d_j(a)}{d_0(a)}E_n^{\frac{1}{2}+\frac{1-j}{m}}+o(1) \underset{n\to\infty}{=}\left\{
\begin{array}{cl}\frac{\left(n-\frac{1}{4}\right)\pi}{d_0(a)},
 &\text{if $\beta=0$,}\\
&\\
\frac{\left(n+\frac{1}{4}\right)\pi}{d_0(a)},
 &\text{if $\beta\not=0$.}
\end{array}\right.
\end{equation}

We then introduce the decomposition
$
E_{n}=E_{n,0}+E_{n,1},
$
where 
\begin{equation}\nonumber
 E_{n,0}=\left\{\begin{array}{lr}
\left(\frac{\left(n-\frac{1}{4}\right)\pi}{d_0(a)}\right)^{\frac{2m}{m+2}},
 &\text{if $\beta=0$,}\\
&\\
\left(\frac{\left(n+\frac{1}{4}\right)\pi}{d_0(a)}\right)^{\frac{2m}{m+2}},
 &\text{if $\beta\not=0$}
\end{array}\right.
\quad\text{and}\quad \frac{E_{n,1}}{E_{n,0}}=o\left(1\right).
\end{equation}
So  we have
\begin{align}
E_{n,0}^{\frac{1}{2}+\frac{1}{m}}
&=E_{n,0}^{\frac{1}{2}+\frac{1}{m}}\left(1+\frac{E_{n,1}}{E_{n,0}}\right)^{\frac{1}{2}+\frac{1}{m}}+\sum_{j=1}^{\lfloor\frac{m}{2}+1\rfloor}\frac{d_j(a)}{d_0(a)}E_{n,0}^{\frac{1}{2}+\frac{1-j}{m}}\left(1+\frac{E_{n,1}}{E_{n,0}}\right)^{\frac{1}{2}+\frac{1-j}{m}}+o(1)\nonumber\\
&=E_{n,0}^{\frac{1}{2}+\frac{1}{m}}\left(1+\sum_{k=1}^{\infty}{\frac{1}{2}+\frac{1}{m}\choose k}\left(\frac{E_{n,1}}{E_{n,0}}\right)^k\right)\nonumber\\
&+\sum_{j=1}^{\lfloor\frac{m}{2}+1\rfloor}\frac{d_j(a)}{d_0(a)}E_{n,0}^{\frac{1}{2}+\frac{1-j}{m}}\left(1+\sum_{k=1}^{\infty}{\frac{1}{2}+\frac{1-j}{m}\choose k}\left(\frac{E_{n,1}}{E_{n,0}}\right)^k\right)+o(1).\nonumber
\end{align}
Thus,
\begin{align}
0&={\frac{1}{2}+\frac{1}{m}\choose 1}\frac{E_{n,1}}{E_{n,0}}+\sum_{k=2}^{\infty}{\frac{1}{2}+\frac{1}{m}\choose k}\left(\frac{E_{n,1}}{E_{n,0}}\right)^k\nonumber\\
&+\sum_{j=1}^{\lfloor\frac{m}{2}+1\rfloor}\frac{d_j(a)}{d_0(a)}E_{n,0}^{-\frac{j}{m}}\left(1+\sum_{k=1}^{\infty}{\frac{1}{2}+\frac{1-j}{m}\choose k}\left(\frac{E_{n,1}}{E_{n,0}}\right)^k\right)+o\left(E_{n,0}^{-\frac{1}{2}-\frac{1}{m}}\right),\nonumber
\end{align}
and hence
\begin{align}
&{\frac{1}{2}+\frac{1}{m}\choose 1}\frac{E_{n,1}}{E_{n,0}}+\sum_{k=2}^{\infty}{\frac{1}{2}+\frac{1}{m}\choose k}\left(\frac{E_{n,1}}{E_{n,0}}\right)^k\nonumber\\
&+\sum_{j=1}^{\lfloor\frac{m}{2}+1\rfloor}\frac{d_j(a)}{d_0(a)}E_{n,0}^{-\frac{j}{m}}\left(\sum_{k=1}^{\infty}{\frac{1}{2}+\frac{1-j}{m}\choose k}\left(\frac{E_{n,1}}{E_{n,0}}\right)^k\right)+o\left(E_{n,0}^{-\frac{1}{2}-\frac{1}{m}}\right)\nonumber\\
&=-\sum_{j=1}^{\lfloor\frac{m}{2}+1\rfloor}\frac{d_j(a)}{d_0(a)}E_{n,0}^{-\frac{j}{m}}.\label{asy_eq4}
\end{align}
Thus, one concludes
$
\frac{E_{n,1}}{E_{n,0}}=E_{n,2}+E_{n,3},
$
where
\begin{equation}\label{ex_eq1}
E_{n,2}=-\frac{2m}{m+2}\frac{d_{1}(a)}{d_{0}(a)}E_{n,0}^{-\frac{1}{m}}\,\,\text{ and }\,\,E_{n,3}=o\left(E_{n,0}^{-\frac{1}{m}}\right).
\end{equation}

 Hence, from \eqref{asy_eq4} and  \eqref{ex_eq1} we have
\begin{align}
&{\frac{1}{2}+\frac{1}{m}\choose 1}\left(E_{n,2}+E_{n,3}\right)+\sum_{{k}=2}^{\infty}{\frac{1}{2}+\frac{1}{m}\choose {k}}\left(E_{n,2}+E_{n,3}\right)^k\nonumber\\
&+\sum_{j=1}^{\lfloor\frac{m}{2}+1\rfloor}\frac{d_j(a)}{d_0(a)}E_{n,0}^{-\frac{j}{m}}\sum_{{k}=1}^{\infty}{\frac{1}{2}+\frac{1-j}{m}\choose {k}}\left(E_{n,2}+E_{n,3}\right)^k+o\left(E_{n,0}^{-\frac{1}{2}-\frac{1}{m}}\right)\nonumber\\
&=-\sum_{j=1}^{\lfloor\frac{m}{2}+1\rfloor}\frac{d_j(a)}{d_0(a)}E_{n,0}^{-\frac{j}{m}}.\label{asy_eq5}
\end{align}
This provides the induction basis.

Next, suppose that 
$\frac{E_{n,1}}{E_{n,0}}=E_{n,2}+E_{n,4}+\cdots+E_{n,2s}+E_{n,2s+1},$
where $E_{n,2s+1}=o\left(E_{n,0}^{-\frac{s}{m}}\right)$ and
$
E_{n,2t}=e_{t}(a)E_{n,0}^{-\frac{t}{m}},\,\, 1\leq t\leq s<\frac{m+2}{2}$ for some  $e_{t}(a)\in\C$.
Then from \eqref{asy_eq4}
\begin{align}
&\sum_{k=1}^{\infty}{\frac{1}{2}+\frac{1}{m}\choose k}\left(E_{n,2}+\cdots+E_{n,2s}+E_{n,2s+1}\right)^k\nonumber\\
&+\sum_{j=1}^{\lfloor\frac{m}{2}+1\rfloor}\frac{d_j(a)}{d_0(a)}E_{n,0}^{-\frac{j}{m}}\sum_{k=1}^{\infty}{\frac{1}{2}+\frac{1-j}{m}\choose k}\left(E_{n,2}+\cdots+E_{n,2s}+E_{n,2s+1}\right)^k+o\left(E_{n,0}^{-\frac{1}{2}-\frac{1}{m}}\right)\nonumber\\
&=-\sum_{j=1}^{\lfloor\frac{m}{2}+1\rfloor}\frac{d_j(a)}{d_0(a)}E_{n,0}^{-\frac{j}{m}}.\nonumber
\end{align}
Hence,
\begin{align}
&\sum_{k=1}^{\infty}{\frac{1}{2}+\frac{1}{m}\choose k}\left(E_{n,2}+\cdots+E_{n,2s}\right)^k\nonumber\\
&+\sum_{j=1}^{\lfloor\frac{m}{2}+1\rfloor}\frac{d_j(a)}{d_0(a)}E_{n,0}^{-\frac{j}{m}}\sum_{k=1}^{\infty}{\frac{1}{2}+\frac{1-j}{m}\choose k}\left(E_{n,2}+\cdots+E_{n,2s}\right)^k\nonumber\\
&=-\sum_{j=1}^{\lfloor\frac{m}{2}+1\rfloor}\frac{d_j(a)}{d_0(a)}E_{n,0}^{-\frac{j}{m}}-{\frac{1}{2}+\frac{1}{m}\choose 1}E_{n,2s+1}+o\left(E_{n,0}^{-\frac{s+1}{m}}\right)+o\left(E_{n,0}^{-\frac{1}{2}-\frac{1}{m}}\right).\label{asy_eq11}
\end{align}
Next, for $1\leq k\geq s+1$
\begin{align}
&\left(E_{n,2}+\cdots+E_{n,2s}+E_{n,2s+1}\right)^k\nonumber\\
&=\left(e_{1}(a)E_{n,0}^{-\frac{1}{m}}+e_{2}(a)E_{n,0}^{-\frac{2}{m}}+\cdots+e_{s}(a)E_{n,0}^{-\frac{s}{m}}+o\left(E_{n,0}^{-\frac{s}{m}}\right)\right)^k\nonumber\\
&=\sum_{k_1=0}^{k}{k\choose k_1}\left(e_{1}(a)E_{n,0}^{-\frac{1}{m}}+e_{2}(a)E_{n,0}^{-\frac{2}{m}}+\cdots+e_{s}(a)E_{n,0}^{-\frac{s}{m}}\right)^{k-k_1}o\left(E_{n,0}^{-\frac{k_1s}{m}}\right)\nonumber\\
&=\left(e_{1}(a)E_{n,0}^{-\frac{1}{m}}+e_{2}(a)E_{n,0}^{-\frac{2}{m}}+\cdots+e_{s}(a)E_{n,0}^{-\frac{s}{m}}\right)^{k}
+o\left(E_{n,0}^{-\frac{s+k-1}{m}}\right)\nonumber\\
&=\sum_{\substack{i_p\geq 0,\,j_p\not=j_q\,\,\text{if}\,\, p\not=q\\i_1+\cdots+i_t=k}}\frac{k!}{i_1!\cdots i_t!}e_{j_1}(a)^{i_1}e_{j_2}(a)^{i_2}\cdots e_{j_t}(a)^{i_t}E_{n,0}^{-\frac{i_1j_1+\cdots+i_tj_t}{m}}+o\left(E_{n,0}^{-\frac{s+k-1}{m}}\right).\nonumber
\end{align}
Also, if $k>s+1$ then $\left(E_{n,2}+\cdots+E_{n,2s}+E_{n,2s+1}\right)^k=o\left(E_{n,0}^{-\frac{s+1}{m}}\right)$.

Then in \eqref{asy_eq11} comparing coefficients of $E_{n,0}^{-\frac{j}{m}}$, $1\leq j\leq s$, we have
\begin{equation}\label{asy_eq12}
-\frac{d_j(a)}{d_0(a)}=\sum_{\substack{|\xi|=k\\ \xi\cdot\eta=j}}{\frac{1}{2}+\frac{1}{m}\choose k}\frac{k!}{\xi!}e(a)^{\xi}+\sum_{r=1}^{j-1}\frac{d_{r}(a)}{d_{0}(a)}\sum_{\substack{|\xi|=k\\ \xi\cdot\eta=j-r}}{\frac{1}{2}+\frac{1-r}{m}\choose k}\frac{k!}{\xi!}e(a)^{\xi},
\end{equation}
where $\eta=(1,2,\dots,m-1)$. 
Moreover, if $\frac{s+1}{m}\leq \frac{1}{2}+\frac{1}{m}$ (i.\ e., $s+1\leq \frac{m+2}{2}$) then there exists some constant $e_{s+1}(a)\in\C$ such that
\begin{equation}\label{asy_eq14}
E_{n,2s+1}=e_{s+1}(a)E_{n,0}^{-\frac{s+1}{m}}+o\left(E_{n,0}^{-\frac{s+1}{m}}\right).
\end{equation}
Now we let $E_{n,2s+2}=e_{s+1}(a)E_{n,0}^{-\frac{s+1}{m}}$ and $E_{n,2s+3}=o\left(E_{n,0}^{-\frac{s+1}{m}}\right)$.

If $s+1> \frac{m+2}{2}$ then $E_{n,0}^{-\frac{s+1}{m}}$ could be smaller than the error term $o\left(E_{n,0}^{-\frac{1}{2}-\frac{1}{m}}\right)$ in \eqref{asy_eq11}, and hence we cannot deduce existence of $e_{s+1}(a)$ like we do in \eqref{asy_eq14}. This completes  induction step and hence proof of Theorem \ref{eigen_asy}.

{\sc email contact:}  kcshin@math.missouri.edu

\begin{thebibliography}{10}

\bibitem{Fedoryuk}
M. V. Fedoryuk.
\newblock {\em Asymptotic Analysis}.
\newblock Springer-Verlag, New York, 1993.


\bibitem{Hille}
E. Hille.
\newblock {\em Lectures on Ordinary Differential Equations}.
\newblock Addison-Wesley, Reading, Massachusetts, 1969.

\bibitem{Shin2}
K. C. Shin.
\newblock Eigenvalues of $\mathcal{PT}$-symmetric oscillators with polynomial potentials.
\newblock {\em Preprint: {\tt math.SP/0407018}}, 23 pages, 2004.

\bibitem{Shin3}
K. C. Shin.
\newblock Schr\"odinger type eigenvalue problems  with polynomial potentials: Asymptotics of  eigenvalues.
\newblock {\em Preprint: {\tt math.SP/0411143}}, 32 pages, 2004.

\bibitem{Sibuya}
Y. Sibuya.
\newblock {\em Global theory of a second order linear ordinary differential equation with a polynomial coefficient}.
\newblock North-Holland Publishing Company, Amsterdam-Oxford, 1975.

\bibitem{TIT}
E. C. Titchmarsh.
\newblock {\em Eigenfunction expansions, Part I}.
\newblock Oxford at the Clarendon Press,  1958.

\bibitem{VOR2}
A. Voros.
\newblock Exact resolution method for general $1$D polynomial Schr\"odinger equation.
\newblock {\em  J. Phys. A: Math. Gen.}, 32:  5993--6007, 1999.

\bibitem{VOR3}
A. Voros.
\newblock   Exercises in exact quantization.
\newblock {\em  J. Phys. A: Math. Gen.}, 33:  7423--7450211--338, 2000.
\end{thebibliography}
\end{document}